\documentclass[a4paper]{amsart}
\usepackage{enumerate}
\usepackage{pb-diagram}
\usepackage{microtype}
\usepackage{amsmath}
\usepackage{amssymb, pb-diagram}
\usepackage[all]{xy}
\usepackage{graphicx} 
\usepackage{mathabx}
\usepackage{chapterbib}
\usepackage{epsfig}
\usepackage{amsfonts}
\usepackage{amssymb}
\usepackage{amsmath}   
\usepackage{amsthm}
\usepackage{color}
\usepackage{latexsym}
\usepackage{hyperref}
\newtheorem{thm}{Theorem}[section]
\newtheorem*{thmi}{Theorem}

\newtheorem{LM}[thm]{Lemma}
\newtheorem{cor}[thm]{Corollary}
 \theoremstyle{definition}

    \newtheorem{rem}[thm]{Remark}

   \DeclareMathOperator {\CAT(0)}{CAT(0)}

  \DeclareMathOperator {\VF}{VF}
  
    \DeclareMathOperator {\Aut}{Aut}

 \DeclareMathOperator {\rank}{rank} 
  
 \DeclareMathOperator {\SL}{SL}

 \DeclareMathOperator {\Ga}{\Gamma} 
    
  \DeclareMathOperator {\Fi}{F}

  \DeclareMathOperator {\FFF}{F}

   \DeclareMathOperator {\Z}{\mathbb{Z}}

 \DeclareMathOperator {\FP}{FP} 
 \DeclareMathOperator {\CAT(0)}{CAT(0)}

    \DeclareMathOperator {\CW}{CW}

\numberwithin{equation}{section}

\begin{document}

\title[Bounding the homological finiteness length]{Bounding the homological finiteness length}

\author{Giovanni Gandini}

\address{School of Mathematics, University of Southampton, Southampton, SO17 1BJ UNITED KINGDOM}
\curraddr{Rheinische Wilhelms-Universit\"{a}t Bonn, Mathematisches Institut, Endenicher Allee 60, 53115 Bonn, Germany}
\email{giovanni.gandini@hausdorff-center.uni-bonn.de}

\date{\today}



\begin{abstract}
We give a criterion for bounding the homological finiteness length of certain ${\scriptstyle \mathbf H}\mathfrak{F}$-groups.  
This is used in two  distinct  contexts. Firstly,  the homological finiteness length of a non-uniform lattice on a locally finite $n$-dimensional contractible $\CW$-complex is less than $n$. In dimension two it solves  a conjecture of Farb, Hruska and Thomas. As another corollary, we obtain an upper bound for the homological finiteness length of arithmetic groups over function fields. This gives an easier proof of a result of Bux and Wortman that solved a long-standing conjecture.  
Secondly, the criterion is applied to integer polynomial points of simple groups over number fields, obtaining bounds established in earlier works of  Bux, Mohammadi and Wortman, as well as new bounds. Moreover, this verifes  a conjecture of Mohammadi and Wortman.
\end{abstract}

\maketitle
\section{introduction}
In 1984, Brown and Geoghegan showed that Thompson's group $\bold{F}$ is of type $\Fi_{\infty}$ giving the first example of a torsion-free group of type $\Fi_{\infty}$ of infinite cohomological dimension. 
 In 1993 Kropholler introduced the class ${\scriptstyle \mathbf H}\mathfrak F$ of \emph{hierarchically decomposable groups} and proved that   torsion-free ${\scriptstyle \mathbf H}\mathfrak F$-groups of type $\FP_{\infty}$ have finite cohomological dimension \cite{MR1246274}. Hence, pathological examples like the group  $\bold{F}$ cannot live inside ${\scriptstyle \mathbf H}\mathfrak F$.  The class ${\scriptstyle \mathbf H}\mathfrak F$ is defined as the smallest class of groups containing the class of finite groups and which contains a group $G$ whenever there is an admissible action of $G$ on a finite-dimensional contractible $\CW$-complex for which all isotropy groups already belong to ${\scriptstyle \mathbf H}\mathfrak F$.  Note that the class ${\scriptstyle \mathbf H}\mathfrak F$ is very large and there are only a few groups which are known to lie outside ${\scriptstyle \mathbf H}\mathfrak F$ \cite{MR1246274, team, gg}.
 
The homological finiteness length $\phi(G)$  of a group $G$ is a generalisation of the concepts of finite generability and finite presentability. 
The main result of  this paper is a  bound for  the homological finiteness length of certain ${\scriptstyle \mathbf H}\mathfrak{F}$-groups: 
\begin{thmi}\label{main} Let $G$ be an ${\scriptstyle \mathbf H}\mathfrak{F}$-group that acts on an $n$-dimensional contractible $\CW$-complex with stabilisers of type $\FP_{\infty}$. If $G$ has no bound on the orders of its finite subgroups or it has infinite Bredon geometric dimension, then $\phi(G)<n$.
\end{thmi}

Suppose now that $X$ is a  locally finite contractible finite-dimensional $\CW$-complex. The topology of uniform convergence  on compact subsets of X turns $\Aut(X)$ into a  locally compact group. 
Every lattice (see Section \ref{S3}) on $X$ lies in  ${\scriptstyle \mathbf H}\mathfrak F$ and we are able to  apply the Theorem to the non-uniform ones. The $S$-arithmetic subgroups of non-commutative almost absolutely simple isotropic algebraic groups over global function fields are important examples of non-uniform lattices. An upper bound for $\phi(\Ga)$ in this  case  is given in \cite{BW}. Our result provides the same bound as a corollary and we are able to prove a slightly strengthened version of a conjecture of  Farb, Hruska, and Thomas \cite[Conjecture 35]{FHT}.

Bux, Mohammadi and Wortman investigated the homological finiteness length of $\SL_{2}(\Z[t, t^{-1}])$, $\SL_{n}(\Z[t])$ and more generally  integer polynomial points of simple groups over number fields \cite{BW, BMW, MW}. 
These groups act on finite-dimensional contractible $\CW$-complexes with arithmetic stabilisers. Note that a group  lies in ${\scriptstyle \mathbf H}\mathfrak F$ whenever it acts admissibly on a finite-dimensional contractible cell-complex with ${\scriptstyle \mathbf H}\mathfrak F$-stabilisers. Hence the groups considered by Bux, Mohammadi and Wortman lie in ${\scriptstyle \mathbf H}\mathfrak F$ and  we can apply our theorem also in this context. In particular we give the same bounds obtained in \cite{BW, BMW, MW} and solve in positive a conjecture of Mohammadi and Wortman \cite{MW}.

\subsection*{Acknowledgements} The author would like to thank Kai--Uwe Bux, Stefan Witzel and Kevin Wortman for the interesting discussions and comments without which this note would not be possible. Moreover, the author is grateful to his Ph.D. supervisor Brita E. A. Nucinkis for her encouragement and advice.

\section{Finiteness properties}
A group $G$ is of type $\FP_{n}$  if the trivial $\Z G$-module $\Z$ admits a resolution of finitely generated projective $\Z G$-modules up to dimension $n$. If $G$ is of type $\FP_{n}$ for every $n\geq 0$, then $G$ is said to be of type $\FP_{\infty}$. 
A group is of type $\FFF_{n}$ if it admits a $K(G, 1)$ with finite $n$-skeleton, and $G$ is of type $\FFF_{\infty}$ if it is of type $\FFF_{n}$ for every $n\geq 0$.
For a group, being finitely generated is equivalent to being of type $\FP_{1}$. A group is finitely presented if and only if it is of type $\FFF_{2}$.   
Bestvina and Brady  show the existence of  non-finitely presented groups of type $\FP_{2}$  \cite{bb-97}.

The \emph{homological finiteness length} of $G$ is defined as $$\phi(G) = \sup\{ m | \,\mbox{$G$ is of type $\FP_{m}$\}}.$$

A group  acts \emph{admissibly} on a $\CW$-complex $X$ if the  set-wise stabiliser of each cell coincides with its point-wise stabiliser.  
A $\CW$-complex $X$ is a $G$-$\CW$-complex if $G$ acts admissibly on $X$.

Kropholler's class  ${\scriptstyle \mathbf H}\mathfrak{X}$ admits a hierarchy defined as follows.
A group $G$ belongs to ${\scriptstyle \mathbf H}_{1}\mathfrak{X}$ if there exists a finite-dimensional contractible $G$-$\CW$-complex $X$ with cell stabilisers in $\mathfrak{X}$.  
The class ${\scriptstyle \mathbf H}_{\alpha}\mathfrak{X} $  is defined by transfinite recursion: 
\begin{itemize}
\item $ {\scriptstyle \mathbf H}_{0}\mathfrak{X}= \mathfrak{X}$; 
\item  if $\alpha$ is a successor ordinal then $ {\scriptstyle \mathbf H}_{\alpha}\mathfrak{X} = {\scriptstyle \mathbf H}_{1} ({\scriptstyle \mathbf H}_{\alpha -1} \mathfrak{X})$; 
\item if $\alpha$ is a limit ordinal then $ {\scriptstyle \mathbf H}_{\alpha}\mathfrak{X} = \bigcup_{\beta < \alpha} {\scriptstyle \mathbf H}_{\beta} \mathfrak{X}$. 
\end{itemize}The operator ${\scriptstyle \mathbf H}$ is defined by ``$G$ belongs to ${\scriptstyle \mathbf H}\mathfrak{X}$ if and only if $G$ belongs to ${\scriptstyle \mathbf H}_{\alpha} \mathfrak{X}$ for some ordinal $\alpha$''.

The next lemma is a well-known criterion for finiteness that follows from \cite[Theorem 2.2]{brownfi}.
\begin{LM}[\cite{brownfi}]\label{brow} Let $G$ be a group that acts on an $n$-dimensional contractible $\CW$-complex with stabilisers of type $\FP_{\infty}$. Then, $G$ is of type $\FP_{\infty}$ if and only if it is of type $\FP_{n}$.
\end{LM}

The next result, due to Kropholler and Mislin, is the main ingredient in the proof of the Theorem.
\begin{thm}[\cite{MR1246274,MR1610595}]\label{krop}
Every ${\scriptstyle \mathbf H}\mathfrak{F}$-group of type $\FP_{\infty}$ has a bound on the orders of its finite subgroups and admits a finite-dimensional classifying space for proper actions.
\end{thm}
The \emph{Bredon geometric dimension} of a group $G$  is the minimal dimension of a classifying space for proper actions $\underline{E}G$. Note that for a torsion-free group the Bredon geometric  dimension coincides with the geometric dimension.
\begin{proof}[Proof (of Theorem)]By an application of Theorem \ref{krop} we obtain that $G$ is not of type $\FP_{\infty}$. Now we apply   Lemma \ref{brow} to conclude that $G$ is not of type $\FP_{n}$.
\end{proof}

\section{Non-uniform lattices on locally finite  contractible cell complexes}\label{S3}

A subgroup $\Gamma$ of a locally compact topological group $G$ with left-invariant Haar measure $\mu$ is a \emph{lattice} if:
\begin{itemize}
\item $\Gamma$ is discrete, and
\item $\mu(\Gamma \backslash G)< \infty$.
\end{itemize}
A lattice $\Ga$ is said to be uniform if $\Ga \backslash G$ is compact.  
 Let $X$ be a locally finite contractible finite-dimensional $\CW$-complex and let $\Aut(X)$ be its full group of  admissible  automorphims.
Note that  every subgroup $G \leq \Aut(X)$ acts admissibly  on the barycentric subdivision of $X$. Moreover, $\Aut(X)$ is locally compact and so it is reasonable to  talk about lattices on locally finite contractible finite-dimensional  $\CW$-complexes. 
\begin{LM}[\cite{treelattices}, Corollary 1.7]\label{clas} Let $X$ be a locally finite  contractible finite dimensional $\CW$-complex with vertex set $V(X)$. If $\Ga$ is a subgroup of $G=\Aut(X)$, then:
\begin{itemize}
\item $\Gamma$ is discrete if and only if the stabiliser $\Ga_{x}$ is finite for each $x \in V (X)$;
\item $\mu(\Ga\backslash G)<\infty$ if and only if  $\sum_{\sigma\in \Ga\backslash X}\frac{1}{|\Ga_{\tilde{\sigma}}|}<\infty$,  where  $\sigma= [\tilde{\sigma}]$.    Moreover,  the Haar measure $\mu$ can be normalised in such  a way that for every discrete $\Ga\leq G$, $\mu(\Ga \backslash G) = \sum_{\sigma\in \Ga\backslash X}\frac{1}{|\Ga_{\tilde{\sigma}}|}$.
  \end{itemize}
\end{LM}

\begin{thm}\label{main} If $\Ga$ is a non-uniform lattice on a locally finite contractible $\CW$-complex of dimension $n$, then   $\phi(\Ga)< n$.
\begin{proof} 
Let $\Ga$ be a non-uniform lattice on a locally finite contractible $\CW$-complex $X$ of dimension $n$. 
By Lemma \ref{clas} $\mu(\Ga\backslash\Aut(X))=\sum_{\sigma\in \Ga\backslash X}\frac{1}{|\Ga_{\tilde{\sigma}}|}$.
Since $\Ga$ is non-uniform, the set $\Ga \backslash X$ is infinite and so  for any $m$ there is some $\sigma \in\Ga \backslash X$ such that $\frac{1}{|\Ga_{\tilde{\sigma}}|}< \frac{1}{m}$. Therefore,  there is no bound on the orders of the stabilisers (which are finite), and so there is no bound on the orders of the finite subgroups of $\Ga$.
Since $\Gamma \in {\scriptstyle \mathbf H}_{1}\mathfrak{F}$ an application of the Theorem gives the result.  
Alternatively this final step can be achieved by noticing that the rational  cohomological dimension of $\Gamma$ is at most $n$ and applying \cite[Proposition]{MR1246274} and Lemma \ref{brow}  or directly  by \cite[Proposition 1]{bound}.
\end{proof} 
\end{thm}
A first immediate application of Theorem \ref{main} is a classical result \cite[Pg. 32]{treelattices}. 
\begin{cor}
If  $X$ is a tree, then every non-uniform lattice in $\Aut(X)$ is not finitely generated.  
More generally, if $\Gamma$ is a   non-uniform lattice on a product of $n$ trees, then $\phi(\Gamma)<n$.
\end{cor}

Suppose now that $X$ is a locally finite $\CAT(0)$ polyhedral complex \cite{FHT}, since every  $\CAT(0)$ space is contractible \cite[Corollary
II. 1.5]{BH} we obtain the following:

\begin{cor} If $\Ga$ is a non-uniform lattice on a locally finite $\CAT(0)$ polyhedral  complex of dimension $n$, then   $\phi(\Ga)< n$.
 \end{cor}
 \begin{rem}
If a group $G$ acts on a $\CW$-complex $X$ with finite stabilisers and $X$  admits a compatible $\CAT(0)$-metric, then $X$ is a model for $\underline{E}G$ \cite[Proposition 3]{bln}.
\end{rem}
\begin{cor} A non-uniform lattice on a locally finite $2$-dimensional contractible $\CW$-complex cannot be finitely presented.
\end{cor}
The above solves a conjecture of Farb, Hruska, and Thomas  \cite[Conjecture 35]{FHT}. 

As a final corollary to Theorem \ref{main} we obtain the main
result of \cite{BW}.
We start by recalling some standard nomenclature. Let $K$ be a global function field, and $S$ be a finite non-empty
set of pairwise inequivalent valuations on $K$.  Let $\mathcal{O}_{S}\leq K$ be the
ring of $S$-integers. Denote a reductive K-group by $\mathbf{G}$. Given a valuation $v$ of $K$, $K_{v}$ is the completion of $K$
with respect to $v$. If $L/K$ is a field extension, the $L$-$\rank$ of $\mathbf{G}$, $\rank_{L} \mathbf{G}$ is the dimension of a maximal $L$-split torus of $\mathbf{G}$. The $K$-group $\mathbf{G}$ is $L$-isotropic if $\rank_{L}\mathbf{G}\neq 0$. As in  \cite{BW}, to any $K$-group $\mathbf{G}$, there is associated a non-negative integer $k(\mathbf{G}, S) = \sum_{v\in S}\rank_{K_{v}} \mathbf{G}$.

 \begin{cor}[Theorem 1.2, \cite{BW}]\label{1.2}
Let  $\bold{H}$ be a connected non-commutative absolutely almost simple $K$-isotropic $K$-group. Then $\phi(\bold{H}(\mathcal{O}_{S})) \leq k(\bold{H}, S) -1$.
\begin{proof}Let $\bold{H}$  be a connected non-commutative absolutely almost simple $K$-isotropic $K$-group.   Let $H$ be  $\prod_{v\in S}\bold{H}(K_{v})$,  there is  a $k(\bold{H}, S)$-dimensional Euclidean building $X$ associated to $H$. $X$  is a locally finite $\CAT(0)$ polyhedral complex.  
The arithmetic group $\bold{H}(\mathcal{O}_{S})$ becomes  a lattice of $H$ via the diagonal embedding since $\bold{H}(\mathcal{O}_{S})$ has finite covolume  by \cite[pg. 41]{harder}. Moreover, $\bold{H}$ is $K$-isotropic if and only if $\bold{H}(\mathcal{O}_{S})$ is non-uniform by \cite[Corollary 2.2.7]{harder}. An application of Theorem \ref{main} completes the proof.
\end{proof}
\end{cor}

\begin{rem} Theorem \ref{main} gives the upper bound  on the homological finiteness length  of arithmetic groups over function fields,  a historical overview can be found in \cite{BW}.
In a recent remarkable paper \cite{bgw} Bux, Gramlich and Witzel   showed that $\phi(\bold{H}(\mathcal{O}_{S})) = k(\bold{H}, S) -1$.  Calculating the homological finiteness length of non-uniform lattices on $\CAT(0)$ polyhedral complexes is an ambitious open problem. We conclude by mentioning that Thomas and Wortman exhibit examples of non-finitely generated non-uniform  lattices on regular right-angled buildings \cite{T-W}. This shows that the upper bound of Theorem \ref{main} is not sharp and in particular, that the Theorem of Bux, Gramlich and Witzel does not hold for all non-uniform lattices on locally finite $\CAT(0)$ polyhedral complexes.
\end{rem}
\section{Actions on non-locally finite buildings}
\begin{cor} The following bounds hold:
\begin{enumerate}
\item The group $\SL_{n}(\Z[t])$ is not of type $\FP_{n-1}$
\item The group $\SL_{n}(\Z[t, t^{-1}])$ is not of type $\FP_{2(n-1)}$
\item Let $K$ be a number field and let $\mathcal{O}_{K}$ be its ring of integers.  Let $\bold{H}$ be a connected, noncommutative, absolutely almost simple algebraic $K$-group whose $K$-rank equals $k$. Then $\bold{H}(\mathcal{O}_{K}[t])$ is not of type  $\FP_k$.
\end{enumerate}
\begin{proof} The group $\SL_{n}(\Z[t])$ acts on a $(n-1)$-dimensional Bruhat-Tits building $X$.  By \cite[Lemma 2]{BMW} the $\SL_{n}(\Z[t])$ stabilisers of cells in $X$ are arithmetic groups. Arithmetic  groups are of type $\VF$ by \cite[Theorem 9.3]{BS73}. By a construction of Serre every group of type $\VF$ has finite Bredon geometric dimension   \cite{serreco} and is of type $\Fi_{\infty}$. In particular, arithmetic groups lie in   ${\scriptstyle \mathbf H}_{1}\mathfrak{F}$ and so $\SL_{n}(\Z[t])\in {\scriptstyle \mathbf H}_{2}\mathfrak{F}$. 
It is clear that $\SL_{n}(\Z[t])$ contains a free-abelian group of infinite countable rank, therefore the Bredon geometric dimension of $\SL_{n}(\Z[t])$ is infinite. An application of Theorem \ref{main} gives  part $(1)$.
It remains to prove part $(2)$ and  part $(3)$. Arguing as in \cite[Section 3]{BWfi} it is possible to show that the  $\SL_{n}(\Z[t, t^1])$ stabilisers of cells in the $2(n-1)$-dimensional Bruhat-Tits building are arithmetic groups as remarked at the end of \cite[Section 4]{BWfi}. In the proof of  \cite[Lemma 2]{MW} the authors show that the  $\bold{H}(\mathcal{O}_{K}[t])$  stabilisers of cells in the $k$-dimensional Bruhat-Tits building are arithmetic groups. In their proof there is no need to assume that the $K$-rank of $\bold{H}$ equals  $2$. Now proceed as in part $(1)$ to obtain the desired bounds.  
\end{proof}
\end{cor}
\begin{rem}Part $(1)$ was originally shown by Bux, Mohammadi and Wortman in \cite{BMW}. Part $(2)$ for $n=2$ follows from a result of Krsti\'c and McCool \cite{KMccool} and an alternative geometric proof is given by Bux and Wortman \cite{BWfi}.
In the case of  $K$-rank equals $2$ part $(3)$ was proved by Mohammadi and Wortman in \cite{MW}. In the same work they conjecture the validity of part $(3)$ for higher ranks.

\end{rem}
\bibliographystyle{amsplain}
\bibliography{math}

\providecommand{\bysame}{\leavevmode\hbox to3em{\hrulefill}\thinspace}
\providecommand{\MR}{\relax\ifhmode\unskip\space\fi MR }
\providecommand{\MRhref}[2]{%
  \href{http://www.ams.org/mathscinet-getitem?mr=#1}{#2}
}
\providecommand{\href}[2]{#2}
\begin{thebibliography}{10}

\bibitem{team}
Goulnara Arzhantseva, Martin~R. Bridson, Tadeusz Januszkiewicz, Ian~J. Leary,
  Ashot Minasyan, and Jacek {\'S}wi{\c{a}}tkowski, \emph{Infinite groups with
  fixed point properties}, Geom. Topol. \textbf{13} (2009), no.~3, 1229--1263.
  \MR{2496045 (2010b:20069)}

\bibitem{treelattices}
Hyman Bass and Alexander Lubotzky, \emph{Tree lattices}, Progress in
  Mathematics, vol. 176, Birkh\"auser Boston Inc., Boston, MA, 2001, With
  appendices by Bass, L. Carbone, Lubotzky, G. Rosenberg and J. Tits.
  \MR{1794898 (2001k:20056)}

\bibitem{bb-97}
Mladen Bestvina and Noel Brady, \emph{Morse theory and finiteness properties of
  groups}, Invent. Math. \textbf{129} (1997), no.~3, 445--470. \MR{MR1465330
  (98i:20039)}

\bibitem{BS73}
A.~Borel and J.-P. Serre, \emph{Corners and arithmetic groups}, Comment. Math.
  Helv. \textbf{48} (1973), 436--491, Avec un appendice: Arrondissement des
  vari{\'e}t{\'e}s {\`a} coins, par A. Douady et L. H{\'e}rault. \MR{0387495
  (52 \#8337)}

\bibitem{bln}
Noel Brady, Ian~J. Leary, and Brita E.~A. Nucinkis, \emph{On algebraic and
  geometric dimensions for groups with torsion}, J. London Math. Soc. (2)
  \textbf{64} (2001), no.~2, 489--500. \MR{MR1853466 (2002h:57007)}

\bibitem{BH}
Martin~R. Bridson and Andr{\'e} Haefliger, \emph{Metric spaces of non-positive
  curvature}, Grundlehren der Mathematischen Wissenschaften [Fundamental
  Principles of Mathematical Sciences], vol. 319, Springer-Verlag, Berlin,
  1999. \MR{1744486 (2000k:53038)}

\bibitem{brownfi}
Kenneth~S. Brown, \emph{Finiteness properties of groups}, Proceedings of the
  {N}orthwestern conference on cohomology of groups ({E}vanston, {I}ll., 1985),
  vol.~44, 1987, pp.~45--75. \MR{MR885095 (88m:20110)}

\bibitem{bgw}
K.-U. {Bux}, R.~{Gramlich}, and S.~{Witzel}, \emph{{Higher finiteness
  properties of reductive arithmetic groups in positive characteristic: the
  rank theorem}}, ArXiv e-prints (2011).

\bibitem{BMW}
Kai-Uwe Bux, Amir Mohammadi, and Kevin Wortman, \emph{{${\rm SL}_n(\Bbb Z[t])$}
  is not {${\rm FP}_{n-1}$}}, Comment. Math. Helv. \textbf{85} (2010), no.~1,
  151--164. \MR{2563684 (2010k:20066)}

\bibitem{BWfi}
Kai-Uwe Bux and Kevin Wortman, \emph{A geometric proof that {${\rm SL}_2({\Bbb
  Z}[t,t^{-1}])$} is not finitely presented}, Algebr. Geom. Topol. \textbf{6}
  (2006), 839--852 (electronic). \MR{2240917 (2007d:20080)}

\bibitem{BW}
\bysame, \emph{Finiteness properties of arithmetic groups over finiteness
  properties of arithmetic groups over function fields}, Inventiones
  Mathematicae \textbf{167} (2007), 355--378.

\bibitem{FHT}
Benson Farb, Chris Hruska, and Anne Thomas, \emph{Problems on automorphism
  groups of nonpositively curved polyhedral complexes and their lattices},
  Geometry, rigidity, and group actions, Chicago Lectures in Math., Univ.
  Chicago Press, Chicago, IL, 2011, pp.~515--560. \MR{2807842}

\bibitem{gg}
G.~{Gandini}, \emph{{Cohomological invariants and the classifying space for
  proper actions}}, To appear in Groups, Geometry and Dynamics.

\bibitem{harder}
G.~Harder, \emph{Minkowskische {R}eduktionstheorie \"uber
  {F}unktionenk\"orpern}, Invent. Math. \textbf{7} (1969), 33--54. \MR{0284441
  (44 \#1667)}

\bibitem{MR1246274}
Peter~H. Kropholler, \emph{On groups of type {$({\rm FP})\sb \infty$}}, J. Pure
  Appl. Algebra \textbf{90} (1993), no.~1, 55--67. \MR{MR1246274 (94j:20051b)}

\bibitem{MR1610595}
Peter~H. Kropholler and Guido Mislin, \emph{Groups acting on finite-dimensional
  spaces with finite stabilizers}, Comment. Math. Helv. \textbf{73} (1998),
  no.~1, 122--136. \MR{MR1610595 (99f:20086)}

\bibitem{KMccool}
Sava Krsti{\'c} and James McCool, \emph{The non-finite presentability of {${\rm
  IA}(F_3)$} and {${\rm GL}_2({\bf Z}[t,t^{-1}])$}}, Invent. Math. \textbf{129}
  (1997), no.~3, 595--606. \MR{1465336 (98h:20053)}

\bibitem{bound}
Ian~J. Leary and Brita E.~A. Nucinkis, \emph{Bounding the orders of finite
  subgroups}, Publ. Mat. \textbf{45} (2001), no.~1, 259--264. \MR{MR1829588
  (2002b:20074)}

\bibitem{MW}
A.~{Mohammadi} and K.~{Wortman}, \emph{{On presentations of integer polynomial
  points of simple groups over number fields}}, ArXiv e-prints (2011).

\bibitem{serreco}
Jean-Pierre Serre, \emph{Cohomologie des groupes discrets}, Prospects in
  mathematics ({P}roc. {S}ympos., {P}rinceton {U}niv., {P}rinceton, {N}.{J}.,
  1970), Princeton Univ. Press, Princeton, N.J., 1971, pp.~77--169. Ann. of
  Math. Studies, No. 70. \MR{MR0385006 (52 \#5876)}

\bibitem{T-W}
Anne Thomas and Kevin Wortman, \emph{Infinite generation of non-cocompact
  lattices on right-angled buildings}, Algebr. Geom. Top. (2011).

\end{thebibliography}
\end{document}